\documentclass[12pt]{amsart}
\usepackage{amsmath}\usepackage{MnSymbol}
\usepackage{amsfonts}\usepackage{verbatim}\usepackage{bbold}
\usepackage{marvosym}\usepackage{ifsym}
\usepackage{dictsym}
\usepackage{wasysym}
\usepackage{amstext}
\usepackage{amsbsy}
\usepackage{amsopn}\usepackage{hyperref}
\usepackage{amsthm}\usepackage{color}

\theoremstyle{definition}

\def\proclaim#1{\vskip0.5em\noindent{\bf #1}\it}
\def\endproclaim{\vskip0.5em\par\noindent\rm}

\def\proclaim#1{\vskip0.5em\noindent{\bf #1}\it}
\def\endproclaim{\vskip0.5em\par\noindent\rm}

\def\demo#1{\vskip0.5em\noindent{\bf #1\ }}

\def\text#1{\mbox{#1}}
\def\flushpar{\par\noindent}

\def\tag#1{\eqno{(#1)}}
\def\mod{\mbox{ mod }}
\newcommand{\mapright}[1]{%
    \smash{\mathop{%
        \hbox to 1cm{\rightarrowfill}
        }
    \limits^{#1}
    }
}
\newcommand{\mapleft}[1]{%
    \smash{\mathop{%
        \hbox to 1cm{\rightarrowfill} $\frak w$ $\D$-nomalized
        }
    \limits_{#1}
    }
}

\def\e{\epsilon}
\def\a{\alpha}
\def\b{\beta}
\def\G{\Gamma}
\def\g{\gamma}
\def\d{\delta}
\def\D{\Delta}
\def\s{\sigma}
\def\Si{\Sigma}
\def\th{\theta}

\def\x{\times}
\def\ov{\overline}
\def\f{\flushpar}
\def\u{\underline}

\def\B{\mathcal B}

\def\({\left(}
\def\){\right)}\def\bdy{\partial}

\def\<{\langle}
\def\>{\rangle}
\def\bul{\smallskip\f$\bullet\ \ \ $}
\def\Smi{\smallskip\f{ \Smiley\ \ \ }}
\def\lfl{\lfloor}\def\rfl{\rfloor}\def\sms{\smallskip\f}\def\Par{\smallskip\f\P}
\def\sms{\smallskip\f}\def\sbul{\f$\bullet\
\ \ $}\def\sms{\smallskip\f}\def\Smi{\smallskip\f{\large\bf\Smiley\ \ \ }}
\def\Lra{\Longrightarrow}
\def\lfl{\lfloor}\def\rfl{\rfloor}\def\st{\text{such that}}

\begin{document}

 \title{On the bounded cohomology of ergodic group actions.}
\author{ Jon. Aaronson and Benjamin Weiss}
\address[Aaronson]{School of Math. Sciences, Tel Aviv University,
69978 Tel Aviv, Israel.}
\email{aaro@tau.ac.il}
\address[Weiss]{Institute of Mathematics
    Hebrew Univ. of Jerusalem,
    Jerusalem 91904, Israel}
\email{weiss@math.huji.ac.il}
\begin{abstract}In this note we show existence of bounded, continuous, transitive cocycles over a
transitive action by homeomorphisms of any finitely generated group on a Polish space, and
bounded, measurable, ergodic cocycles over any ergodic, probability preserving action of $\Bbb Z^d$.
\end{abstract}
\thanks{\copyright 2017. Aaronson{\tiny$^\prime$}s research was partially supported by ISF grant No. 1289/17.}
\maketitle\markboth{Bounded cohomology}{Aaronson and Weiss}
\section*{\S0 Introduction}
\

\subsection*{Cocycles and skew product actions}
\

Let $\G$ be a countable group and let $X$ be a {\tt space}. In the sequel, $X$ will represent either a Polish, metric space $X=(X,d)$
or a standard probability space $X=(X,\B,m)$. A $\G$-{\it action} on $X$ is an homomorphism
$T:\G\to\text{\tt Aut}\,(X)$. In the topological case, $\text{\tt Aut}\,(X)=\text{\tt Homeo}\,(X)$ in  and
$\text{\tt Aut}\,(X)=\text{\tt PPT}\,(X,\B,m)$ the group of probability preserving transformations of $(X,\B,m)$ in
the probabilistic case.
\

Let $\Bbb G$ be an Abelian topological group equipped with a norm $\|\cdot\|_\Bbb G$.
To define a $\G$-{\tt skew product action} on $X\x\Bbb G$, we need a $T$-{\it cocycle}, that is a
 function $F:\G\x X\to\Bbb G$ satisfying
 \begin{align*}
 \tag{\dsliterary}F(nk,x)=F(k,x)+F(n,T_kx)\ \ (n,\ k\in\G).
 \end{align*}
 The cocycle $F:\G\x X\to\Bbb G$ is assumed to be
continuous in the topological case and measurable in the probabilistic case.
\

 The $F$-{\it skew product} transformations are then defined on $X\x\Bbb G$ by
          $$T^{(F)}_n(x,z):=(T_n(x),z+F(n,x))\ \ (n\in\G).$$
          The assumptions and (\dsliterary) ensure that $T^{(F)}:\G\x (X\x\Bbb G)\to\text{\tt Aut}\,(X\x\Bbb G)$
          is a $\G$-action, called the {\it skew product action}.
          \

          We'll also consider cocycles which are {\it bounded} in the sense that
 $\sup_{x\in X}\|F(\g,x)\|_\Bbb G<\infty\ \forall\ \g\in\G$.
 \

          In case $\G=\Bbb Z$, it is easy to exhibit cocycles. Let
$\varphi:X\to\Bbb G$ and define $F=F^{(\varphi)}:\Bbb Z\x X\to\Bbb G$  by
 $$F(n,x)=\begin{cases} &\ \sum_{k=0}^{n-1}\varphi(T^kx)\ \ \ \ \ \ \ \ n\ge 1,\\ &
0\ \ \ \ \ \ \ \ n=0,\\ & -\sum_{k=1}^{|n|}\varphi(T^{-k}x)\ \ \ \ \ \ \ \ n\le -1.
          \end{cases}$$

          This is a cocycle 
          and indeed any cocycle is of this form. We sometimes write $T^{(F)}_n=T_\varphi^n$
     where     $T_\varphi:X\x\Bbb G\to X\x\Bbb G$ is the {\it skew product transformation}
defined by $T_\varphi(x,z):=(T(x),z+\varphi(x))$.
\

          \

Construction of cocycles for the actions of multidimensional groups (e.g. $\G=\Bbb Z^2$) is more difficult.
\

Note that a constant cocycle for  an action $T:\G\to\text{\tt Hom}\,(X)$ is given by an homomorphism
$h:\G\to\Bbb G$ an in this case, the skew product action $T^{(h)}:\G\to\text{\tt Hom}\,(X\x\Bbb G)$ is given
by the (direct) product action $T\x h$ where $(T\x h)_\g(x,y):=(T_\g(x),y+h(\g)).$

\

The simplest $\Bbb G$-valued, non constant $T$-cocycles for  an action $T:\G\to\text{\tt Hom}\,(X)$ are given by a
{\it coboundary}, that is, a function $h:\G\x X\to\Bbb G$ defined by
\begin{align*}
h(n,x)=c(x)-c(T_nx).
\end{align*} where $c:X\to\Bbb G$ (the {\it transfer function}) is measurable or continuous in the probabilistic and topological cases respectively.
It is not hard to see that a coboundary is a cocycle.
\

For full shifts of $\Bbb Z^d$, the only H\"older continuous $\Bbb R^\kappa$-valued cocycles are sums of a
coboundary and a constant cocycle (homomorphism). See \S4.
The dynamical properties of such cocycles are somewhat limited (see \S4). However for certain infinitely generated groups,
constant cocycles can have robust dynamics (see \S3).
        \subsection*{Topological cocycles}
          \

         The first topologically transitive, topological skew product $\Bbb Z$-actions on
          $\Bbb T\x\Bbb R$ were constructed  in \cite{Sni30} and \cite{Bes37}.
          \

          We prove:
 \

 \proclaim{Theorem 1}\ \ Let $X$ be a perfect, Polish space, $\G$ be a countable, finitely generated group and let
 $(\Bbb B,\|\cdot\|_\Bbb B)$ be a separable Banach space.
 \

 If $T:\G\to\text{\tt Homeo}\,(X)$ is a topologically transitive action, then there is a bounded,  continuous cocycle
 $h:\G\x X\to\Bbb B$ so that the skew product action $(X\x\Bbb B,T_h)$ is topologically transitive.\endproclaim

 Theorem 1 was established for $\Bbb Z$-actions in \cite{Sidorov}. Note that the skew product action of a coboundary cannot be transitive.
 The cocycle in Theorem 1 is constructed as a limit of coboundaries. In \S3, we consider versions of theorem 1 for certain infinitely generated
 groups. 

\subsection*{Measurable cocycles}
\

As shown in \cite{Herman} $\&$ \cite{AW-Herman}, for $(X,\B,m,T)$ an ergodic, probability preserving transformation and $\Bbb G$ a locally compact, Polish, amenable group, there is
bounded, measurable function $F:X\to\Bbb G$ so that the skew product  $(X\x\Bbb G,m\x m_\Bbb G,T_F)$ is ergodic.

 \

Now let $S:\mathbb{\Gamma}\x Y\to\text{\tt Aut}\,(Y)$ be
an ergodic action of the countable, amenable group $\mathbb{\Gamma}$ on the probability space $Y$, there is a   cocycle   $G:\Bbb Z^k\x Y\to\Bbb G$ so that $S^{(G)}$ is ergodic.
 \

 By \cite{C-F-W}, the actions $T:\mathbb{Z}\x X\to\text{\tt Aut}\,(X)\ \&\ S:\mathbb{\Gamma}\x Y\to\text{\tt Aut}\,(Y)$
 are orbit equivalent. The orbit equivalence transports $F$ to an  cocycle
 $G:\mathbb{\Gamma}\x Y\to\Bbb G$  for $S$ so that the actions $T^{(F)}\ \&\ S^{(G)}$ are orbit equivalent.
 Thus $S^{(G)}$ is ergodic.
 \

 However,  regularity properties of $F$ (e.g. boundedness) need not pass to $G$.
 \

We prove:
 \

  \proclaim{Theorem 2}\ \ Let $d,\ D\ge 1$, let $\Bbb G\le\Bbb R^D$ be a closed subgroup of full dimension  and let $\mathcal S\subset\Bbb G$ be finite, symmetric and generating in the sense
 that $\overline{\<\mathcal S\>}=\Bbb G$.
  \

Let  $T$ be an ergodic $\Bbb Z^d$-action of the standard non-atomic probability space  $(X,\B,m)$.
There is a bounded, measurable cocycle $F:\Bbb Z^d\x X\to\Bbb G$ with ergodic skew product action $T^{(F)}$ so that  $F(e_k,\cdot)\in\mathcal S\cup\{0\}\ \forall\   1\le k\le d$.\endproclaim
 Here and throughout, $e^{(d)}_k\in\Bbb R^d\ \ (1\le k\le d)$ with
 $(e^{(d)}_k)_j=1_{[j=k]}$. In other words, $\{e^{(d)}_k:\ 1\le k\le d\}$ is the usual orthonormal basis  for
 $\Bbb R^{d}$. When the dimension $d$ is unambiguous, we suppress it and write $e^{(d)}_k=e_k$.

The proof of theorem 1 is given in \S2 and that of theorem 2 is given in \S3. In both proofs the advertised cocycles are limits of
coboundaries satisfying finite {\tt essential value conditions}.
\

\section*{\S1 Proof of theorem 1}

Suppose that $(X,T)$ is a free action of the countable $\G $ by homeomorphisms on a complete, separable, perfect metric space $(X,d)$ and that $(\Bbb B,\|\cdot\|_\Bbb B)$
is a separable Banach space.
\

\subsection*{Bounded,   continuous cocycles}
\

A function $F:X\to\Bbb B$ is  {\it bounded} if $\sup_{x\in X}\|F(x)\|_\Bbb B<\infty$.
\

We denote the collection of {\it bounded,  continuous} (abbr. {\tt BC}) $\Bbb B$-valued functions by $C_B(X,\Bbb B)$. It is a Banach space when equipped with the supremum norm
$$\|F\|_{\sup}:=\sup_{x\in X}\|F(x)\|_\Bbb B.$$
  We call the cocycle $h:\G \x X\to\Bbb B$  {\tt BC} if for each $\g\in\G,\ x\mapsto h(\g,x)$ is a
  {\tt BC} function $X\to\Bbb B$.
  \

Denote the collection of {\tt BC}, $\Bbb B$-valued   cocycles for $T$ by
$\text{\tt Coc}\,(X,T,\Bbb B)$. Fixing a finite, symmetric set  $\Si$ of generators  for $\G$, we may define
$$\|h\|_{\text{\tt Coc}}:=\max_{\g\in\Si}\sup_{x\in X}\,\|h(\g,x)\|_\Bbb B\ \ \ (h\in \text{\tt Coc}\,(X,T,\Bbb B)). $$
It is not hard to see that  $(\text{\tt Coc}\,(X,T,\Bbb B),\|\cdot\|_{\text{\tt Coc}})$ is a Banach space.

\subsection*{Word metric on $\G$}
\

Define the $\Si$-{\it word length norm} on $\G$ by
$$\|N\|_\Si :=\min\,\{k\ge 1:\ \exists\ \s_1,\s_2,\dots,\s_k\in\Si,\ N=\s_1\s_2\dots\s_k\}.$$
It follows that $\|NN'\|_\Si \le \|N\|_\Si \,+\ \|N'\|_\Si $ and, since $\Si$ is symmetric, $\|N^{-1}\|_\Si =\|N\|_\Si $.
\

The left invariant $\Si$-{\it word metric} on $\G$ is the  distance
$$\rho_\Si (k,\ell):=\|k^{-1}\ell\|_\Si .$$

It is easy to see that for $N\in\G\ \&\ h\in \text{\tt Coc}\,(X,T,\Bbb B)$,
$$|h(N,x)|\le \|h\|_{\text{\tt Coc}}\|N\|_\Si .$$

\subsection*{Coboundaries}
The cocycle $h\in \text{\tt Coc}\,(X,T,\Bbb B)$ is a {\it   {\tt BC}  coboundary} if  for some $c\in C_B(X,\Bbb B)$, (the {\it transfer function})
\begin{align*}\tag{\tt cob} h(n,x)=c(x)-c(T_nx).
\end{align*}
In this case we denote $h:=\nabla c$.

We denote the collection of {\tt BC} coboundaries by
$$\bdy C_B(X,T,\Bbb B):=\{\nabla g:\ g\in C_B(X,T,\Bbb B)\}.$$
\

It is easy to see that $\overline{\bdy  C_B(X,\Bbb B)}\subset\text{\tt Coc}\,(X,T,\Bbb B)$. \ \
The cocycle advertised in theorem 1 will appear in $\overline{\bdy  C_B(X,\Bbb B)}$.

\subsection*{Topological essential values at a transitive point}
\

Let $h\in \text{\tt Coc}\,(X,T,\Bbb B)$, fix a $T$-transitive point $x_0\in X$ and set
$$\frak T_{x_0}:=\{\text{ open nbds of}\  x_0\}.$$
\

\ \ For $U\in\frak T_{x_0},\ t\in\Bbb B,\ n\in\G\ \&\ \e>0$ we'll need the
``{\tt essential value conditions}''
\begin{align*}&\tag*{$\text{\tt EVC}(x_0,U,t,\e,n):$}T_nx_0\in U\ \&\ \|h(n,x_0)-t\|_\Bbb B <\e;\\ &
\tag*{$\text{\tt EVC}(x_0,U,t,0,n):$}T_nx_0\in U\ \&\ h(n,x_0)=t.
\end{align*}

\

Let
$$E(h,x_0)=\{r\in \Bbb B:\ \forall \e>0\ \&\ U\in\frak T_{x_0},\ \exists n\in\G\ \st\
\text{\tt EVC}(x_0,U,t,\e,n)\ \text{holds}\}.$$
\

\proclaim{Proposition 1}
\

\ $E(h,x_0)$ is  a closed subsemigroup of $\Bbb B$ and if $E(h,x_0)=\Bbb B$, then
\begin{align*}
 \tag{\Football}\overline{\{R_n(x_0,0):\ n\in\G\}}=\ X \x \Bbb B.
\end{align*}
where $R$ is the $h$-skew product action.

\endproclaim
\demo{Proof}\ \ Let $s,\ t\in E(h,x_0),\ U\in \frak T_{x_0}\ \&\ \e>0$. We'll show that $\exists\ K\in\G$ so that
\begin{align*}
 \tag*{$\bigstar$}T_Kx_0\in U\ \&\ \|h(K,x_0)-(s+t)\|_\Bbb B<\e.
\end{align*}

By definition, $\exists\ n\in\G$ so that $T_nx_0\in U\ \&\ \|h(n,x_0)-s\|_\Bbb B<\tfrac{\e}2$.
\

Thus, $V:=U\cap T^{-1}_{n}U\cap [\|h(n,\cdot)-s\|_\Bbb B<\tfrac{\e}2]\in\frak T_{x_0}\ \&$ by definition
$\exists\ N\in\G\ \st\ T_{N}x_0\in V\ \&\ \|h(N,x_0)-t\|_\Bbb B<\tfrac{\e}2$.
\

Now $T_Nx_0\in V\subset T^{-1}_{n}U\cap [\|h(n,\cdot)-s\|_\Bbb B<\tfrac{\e}2]$ so
$$T_{nN}x_0=T_n\circ T_N(x_0)\in U\ \&\ \|h(n,T_Nx_0)-s\|_\Bbb B<\tfrac{\e}2.$$
\

It follows that
$$\|h(nN,x_0)-(s+t)\|_\Bbb B\le \|h(N,x_0)-t\|_\Bbb B+\|h(n,T_Nx_0)-s\|_\Bbb B<\e.\ \ \CheckedBox\bigstar$$
It is immediate that $E(h,x_0)$ is closed and that  $$(x_0,t)\in \overline{\{T_h^n(x_0,0):\ n=1,2,...\}}\ \forall\ t\in E(h,x_0).$$
Using this and $T_h$-invariance, we see (\Football) when  $E(h,x_0)=\Bbb B$.\ \ \Checkedbox

\

\proclaim{Lemma 2}
\

Let $V\in\frak T_{x_0},\ s\in\Bbb B,\ \D>0\  \&$ let $F\in \bdy\text{\tt Coc}\,(X,T,\Bbb B)$, then
\f$\exists\ h\in\bdy \text{\tt Coc}\,(X,T,\Bbb B),\ \|h\|_{\text{\tt Coc}}<\D\ \&\ n\in\G$ so that  $F+h$ satisfies
$\text{\tt EVC}(x_0,V,s,\D,n)$.
\endproclaim
\demo{Proof}\ \
\

Suppose that $F_g=c-c\circ T_g$ where $c\in C(X)$ and find $W\in\frak T_{x_0},\ W\subset V$ so that
$$|c(x)-c(y)|<\frac{\D}3\ \forall\ x,\ y\in W.$$
Fix $N\in\G$ so that
$$\frac{\|s\|_\Bbb B}{\|N\|_\Si}<\D\ \ \&\ \ T_N(x_0)\in W.$$
 Next, fix $\mathcal E>0$ so that the sets
 $$\{T_kB(x_0,2\mathcal E):\ \|k\|\le 3\|N\|\}$$ are disjoint. Here, for $x\in X,\ r>0$,
 $B(x,r):=\{y\in X:\ d(y,x)\le r\}$ is the closed $r$-ball centered at $x$.
 \

 Let
 \begin{align*}
   & G(x):=\left(1-\frac{d(B(x_0,\mathcal E),x)}{\mathcal E}\right)_+\ \ \ \ (x\in X);\\ &
   a(k):=\left(1-\frac{\rho_\Si (N,k)}{\|N\|}\right)_+\ \ \ \ \ (k\in\G);
 \end{align*}
 where $d(B,x):=\inf_{y\in B}d(y,x)$ and $a_+:=\max\{a,0\}$ for $a\in\Bbb R$.
It follows that $$G\in  C(X,[0,1]),\ G|_{B(x_0,\mathcal E)}\equiv 1\ \&\ G|_{B(x_0,2\mathcal E)^c}\equiv 0$$
and that
$$|a(k)-a(k\g)|\le \frac1{\|N\|}\ \forall\ \g\in\Si\ \&\ k\in\G.$$

\

Now define $g:X\to\Bbb B,\ g\ge 0$ by
$$g(x)=\begin{cases} &\ sa(k)G(T^{-1}_{k}x)\ \ \ \ \ \ \ \ x\in T_kB(x_0,2\mathcal E),\ \rho_\Si (N,k)\le \|N\|,\\ &
  0\ \ \ \ \ \ \ \ x\in X\setminus\bigcup_{k\in\G,\ \rho_\Si (N,k)\le \|N\|}T_kB(x_0,2\mathcal E). \end{cases}$$

It follows that $g\in  C_B(X,\Bbb B)$.
\

Define  $h(n,x):=g(T_{n}x)-g(x)$. It follows from the above that
\begin{align*}
 \|h(\g,x)\|_\Bbb B\le \frac{\|s\|_\Bbb B}{\|N\|_\Si}<\frac{\D}3\ \forall\ \g\in\Si.
\end{align*}

\

To see that $h$ satisfies $\text{\tt EVC}(x_0,W,s,0,N)$,
we have that  $T_N(x_0)\in W$ and
$$h(N,x_0)=g(T_{N}x_0)-g(x_0)=sa(N)-sa(0)=s.$$
Lastly,
\begin{align*}|F(N,x_0)+h(N,x_0)-s|\le |F(N,x_0)|=|c(x_0)-c(T_N(x_0)|<\frac{\D}3
\end{align*}\ \ $\because\ x_0,\ T_Nx_0\in W$. Thus $F+h$ satisfies $\text{\tt EVC}(x_0,V,s,\D,N)$.\ \Checkedbox

\subsection*{Categorical construction}
\

Let $x_0\in X$ be a $T$-transitive point.  For each $V\in\frak T_{x_0},\ s\in\Bbb B,\ n\in\G\  \&\  \D>0$, let
\begin{align*}
 E(&x_0,V,s,\D):=\\ &\{h\in\bdy\,\text{\tt Coc}\,(X,T,\Bbb B):\ \exists\ n\in\G\ \st\ h\ \text{satisfies\
\ {\tt EVC}}\,(x_0,V,s,\D,n)\}.
\end{align*}

It follows from lemma 2 that $E(x_0,V,s,\D)$ is open and dense in $\ov{\bdy\,\text{\tt Coc}\,(X,T,\Bbb B)}$.

Fix
\bul a decreasing sequence $(U_n:\ n\ge 1),\ U_n\in\frak T_{x_0}$ so that $\forall\ V\in\frak T_{x_0}\ \exists\ N_V$ so that
$U_n\subset V\ \forall\ n\ge N_V$;
\bul $(t_n:\ n\ge 1) \in \Bbb B^\Bbb N$ dense, and taking each value i.o.;
and
\bul $\eta_n>\eta_{n+1}\ \searrow\ 0$.
\

By Baire's theorem,
$$E:=\bigcap_{k=1}^\infty E(x_0,U_k,t_k,\eta_k)$$ is residual in $\ov{\bdy\,\text{\tt Coc}\,(X,T,\Bbb B)}$.
\

By proposition 1, for each $h\in E$ we have that
$E(h,x_0)=\Bbb B\ \&$ hence
\begin{align*}
 \tag{\Football}\overline{\{R_n(x_0,0):\ n\in\G\}}=\ X \x \Bbb B.
\end{align*}
This proves theorem 1.\ \ \Checkedbox

We note that sequential constructions are also possible. See \S3.
\section*{\S2 Proof of theorem 2}
\

We adapt here from \S3 of \cite{A-L-V} the {\it essential value conditions}
or {\tt EVC}'s, which give countably many conditions for the ergodicity of the
the  skew product action  $T^{(\varphi)}:\Bbb Z^d\to\text{\tt MPT}\,(X\x\Bbb G,\B(X\x\Bbb G),m\x m_\Bbb G)$.
\par These are best understood in terms
of  cocycles with respect to the {\tt orbit equivalence relation} of $T$  and its {\tt groupoid}  
as in \cite{F-M}, (see below).
\subsection*{Orbit cocycles}

\par The {\it orbit equivalence relation} generated by the free $\Bbb Z^d$-action  $T$ is
$$\mathcal R=\mathcal R_T:=\{(x,T_nx):x\in X,\ n\in\Bbb Z^d\}.$$

\

An $\mathcal R$-{\it cocycle } is a measurable function $\tilde\varphi:\mathcal R\to
 \Bbb G$ such that if $(x,y),\ (y,z)\in\ \mathcal R$, then
$$\tilde\varphi(x,z)=\tilde\varphi(x,y)+\tilde\varphi(y,z).$$
\sms The $\mathcal R$-cocycle $\tilde\varphi:\mathcal R\to
\Bbb G$ corresponds to a $T$-cocycle $\varphi:\Bbb Z^d\x X\to \Bbb G$ via
 $$\varphi(n,x):=\tilde\varphi(x,T_nx).$$

\subsection*{Groupoid}
\

\par A {\it partial probability preserving transformation} of $X$ is a pair $(R,A)$ where $A\in\B$ and $R:A\to RA$ is measurable, invertible and
$m|_{RA}\circ R^{-1}=m|_A$. The set $A$ is called the {\it domain } of $(R,A)$.
We'll sometimes abuse this notation by writing $R=(R,A)$ and $A=\mathcal D(R)$.
Similarly, the {\it image} of $(R,A)$ is the set $\Im(R)=RA$.
\

An $\mathcal R$-{\it holonomy} is a partial probability preserving transformation $R$ of $X$ with the additional property that
$$(x,R(x))\in\mathcal R\ \forall\ x\in \mathcal D(R).$$
The {\it groupoid of $\mathcal R$} (or of $T$) is the collection
$$[[\mathcal R]]=[[T]]:=\{\text{$\mathcal R$-holonomies}\}.$$
The {\it full group of $\mathcal R$} is
$$[T]:=\{R\in [[\mathcal R]]:\ \mathcal D(R)=\Im(R)=X\ \mod\ m\}.$$
For $R\in [[T]]$, the function $x\mapsto\varphi(R,x)\ \ (\mathcal D(R)\to \Bbb G)$ is defined by
$$\varphi(R,x)=\tilde\varphi(x,Rx).$$
The cocycle property ensures that $\varphi(R\circ S,x)=\varphi(S,x)+\varphi(R,Sx)$ on
$\mathcal D(R\circ S)=\mathcal D(S)\cap S^{-1}\mathcal D(R)$ and for $\varphi_{T_n}(x)=\varphi(n,x)$.
\par An $\mathcal R$-holonomy can be thought of as a {\tt random power of $T$}.
For $A\in\B(X)$ and $\chi:A\to\Bbb Z^d$ measurable, define $T^{(\chi)}:A\to X$ by
$T^{(\chi)}(x):=T_{\chi(x)}x$.  Any $\mathcal R$-holonomy is  of this form (but a "random power" need not be a
partial probability preserving transformation).

\

\par The orbit equivalence relation of the skew product action  $T^{(\varphi)}$ is given by
$$\mathcal R_{T^{(\varphi)}}=\{((x,y),(x',y'))\in (X\x\Bbb G)^2:\ (x,x')\in\mathcal R_T\ \&\ y'=y+\tilde\varphi(x,x')\}.$$

\subsection*{Essential value conditions} Let  $A\in\B,\ U$ a subset of $\Bbb G$, and $c>0$.
We say that the measurable cocycle $\varphi:X\to \Bbb G$
satisfies {\tt EVC}$_T(U,c,A )$ if
 $\exists\ R\in [[T]]$
such that
$$\mathcal D(R),\ \Im(R)\subset A,\ \varphi_R\in U\text{ on }\mathcal D(R),\ m(\mathcal D(R)))>c m(A).$$
This is in honor of the collection of {\it essential values} introduced in \cite{S}
\begin{align*}
 &E(T,\varphi):=\\ &\{a\in\Bbb G:\ \forall\ A\in\B_+,\ U\in\frak T_a,\ \exists\ n\in\Bbb Z^d,\ m(A\cap T_n^{-1}A\cap [\varphi(n,\cdot)\in U])>0\}
\end{align*}
where $\frak T_a:=\{ U\ni\ a\ \text{open in}\ \Bbb G\}$.
\

It is shown in \cite{S} that $E(T,\varphi)$ is a closed subgroup of $\Bbb G$ and that $T^{(\varphi)}$ is ergodic iff
$E(T,\varphi)=\Bbb G$. The following is a standard consequence of this. See \cite{A-L-V} or \cite{AW-Herman}.
\proclaim{Ergodicity Proposition}
\par The skew
product action
$T^{(\varphi)}$
is
ergodic with respect to the product measure $m\x m_{ \Bbb G}$ iff $\exists$
\bul a countable
 base $\mathcal U$ for $\frak T_0$;
\bul a countable, dense collection $\mathcal A\subset\B$;
\bul a countable  collection $\mathcal S\subset\Bbb G$ so that $\overline{\<\mathcal S\>}=\Bbb G$
\par and a number $0<c<1$ such that
\f $\varphi$ satisfies {\tt EVC}$_T(\sigma+U,c,A )\ \forall\ A\in \mathcal A,\ \sigma\in\mathcal S,\ U\in\mathcal U$.
\endproclaim

\par Essential value conditions are impervious to small changes. The following is a standard modification of lemma 3.5 of \cite{A-L-V}.
\proclaim{Stability Lemma}
\par If
$\psi:\Bbb Z^d\x X\to \Bbb G$ is a  measurable cocycle
satisfying {\tt EVC}$_T(U,c,A)$ where $A\in\B,\ c >0,\ U\subset \Bbb G$;
then $\exists\ \d>0$ such that
if $\varphi:\Bbb Z^d\x X\to \Bbb G$ is another measurable cocycle, and
$$m([\varphi(e_k,\cdot)\ne \psi(e_k,\cdot)])<\d\ \forall\ 1\le k\le d,$$
then $\varphi$ satisfies
{\tt EVC}$_T(U,c,A)$.
\endproclaim

\subsection*{Coboundaries}
\

The advertised cocycle is constructed as a limit of coboundaries,   a {\it coboundary} being a cocycle
$\psi:\Bbb Z^d\x X\to\ \Bbb G$ of form $\psi(n,x)=F(x)-F(T_nx)$ where $F:X\to\Bbb G$ is a measurable function called the
{\it transfer function}. We'll denote the coboundary with transfer function $F$ by
$$\nabla F(n,x):=F(x)-F(T_nx).$$

\subsection*{Discrete distance}
\

Let $\|\cdot\|=\|\cdot\|_1$ on $\Bbb R^d$. We'll consider $\Bbb Z^d$ as a discrete metric space and write
\begin{align*}&\text{For}\ k\in\Bbb Z^d,\ R>0,\ \Si(k,R):=\{j\in\Bbb Z^d:\ \|j-k\|\le R\}\ \&\ \Si_R:=\Si(0,R).\\ &
\text{For}\ Q\subset\Bbb Z^d,\ R>0,\ \Si(Q,R):=\bigcup_{k\in Q}\Si(k,R),\\ &\,\,\,\,\,\,\,\,\ Q^o:=\{k\in Q:\ \Si(k,1)\subset Q\}\ \&\ \bdy Q:=Q\setminus Q^o.
\end{align*}

\subsection*{Rokhlin towers}
\

A {\it Rokhlin tower} is a collection $\mathcal T=\mathcal T_{N,B}=\{T_kB:\ k\in\Si_N\}$ where $B\in\B$ is such that these sets are disjoint.
\

The {\it breadth} of $\mathcal T=\mathcal T_{N,B}$ is $N_\mathcal T=N$, the {\it base} is $B_\mathcal T=B$ $\&$  the {\it error}
of the Rokhlin tower is $\e_\mathcal T:=m(X\setminus\bigcupdot_{k\in\Si_N}T_kB)$.
\

The {\tt Rokhlin Lemma} for $\Bbb Z^d$ actions as in \cite{Conze} and \cite{K-W}, (see also \cite{O-W})
says that
\f{\sl any free, ergodic, probability preserving $\Bbb Z^d$ action has Rokhlin towers
of any breadth  and   error (in $(0,1)$).}
\subsection*{Castles}
A {\it castle} is an array of disjoint Rokhlin towers with the same breadth.
A castle may be derived from a Rokhlin tower by partitioning its base.
\

The {\it interior} of $\mathcal T_{N,B}$ is
$$\mathcal T_{N,B}^o:=\{T_kB:\ k\in\Si_N,\ k\pm e_i\in\Si_N,\ i=1,2,\dots,d\}=\mathcal T_{N-1,B}$$
and the
{\it boundary} of $\mathcal T_{N,B}$ is
$$\bdy\mathcal T_{N,B}:=\mathcal T_{N,B}\setminus\mathcal T_{N,B}^o=\{T_kB:\ \|k\|=N\}.$$
\

Let $Q\subset\Si_N$. We'll write
$$T_QB:=\bigcupdot_{k\in Q}T_kB.$$
\subsection*{Purifications}
\

Given a Rokhlin tower $\mathcal T_{N,B}$ and a partition $\a\subset\B$, the $\a$-{\it purification} of $B$ is
is the partition
$$\b:=\{B_{\u a}:=B\cap\bigcap_{\|k\|\le N}T_k^{-1}a_k:\  \u a\in \a^{\Si_N}\}$$
of $B$ and  the $\a$-{\it purification of $\mathcal T_{N,B}$} is the corresponding castle
$$\mathcal P=\mathcal P_{\mathcal T,\a}:=\{T_kb:\ k\in\Si_N,\ b\in\b\}.$$
\

For the rest of this paper, we fix a finite, symmetric generator set $\mathcal S\subset\Bbb G\setminus\{0\}$.
 \subsection*{Step functions}
\

Let  $\mathcal P=\{T_kb:\ k\in\Si_N,\ b\in\b\}$  be a purification of the Rokhlin tower $\mathcal T=\mathcal T_{N,B}$. A $\mathcal P$-{\it step  function} $F:X\to\Bbb G$
is one of form
$$F=\sum_{k\in\Si_N,\ b\in \b}a_{k,b}1_{T_kb}.$$
It is called $\mathcal T$-{\it internal} if $F|_{T_{\bdy\Si_N}B}\equiv 0$ and  
$\mathcal S$-{\it incremental}
\begin{align*}\nabla F(e_i,\cdot)\in\mathcal S\cup\{0\}\ \ (i=1,2,\dots,d).
\end{align*}

\proclaim{Inductive lemma}
\

Let  $\mathcal P_0$  be a purification of the Rokhlin tower $\mathcal T_0$ and let $F_0:X\to\Bbb G$ be a  $\mathcal T_0$-internal, $\mathcal S$-incremental
$\mathcal P_0$-{ step  function}.
\

Fix $\e>0,\ 0<r<\tfrac1{2^{d+3}},\ \sigma\in\mathcal S$ and $A\in\B_+$. There exist
\sbul a Rokhlin tower $\mathcal T$ with $\e_\mathcal T<\e$ and a purification $\mathcal P$, and
\sbul a $\mathcal T$-internal, $\mathcal S$-incremental  $\mathcal P$-{ step  function} $F:X\to\Bbb G$ so that
\begin{align*}\tag{a}
 \mu([\nabla F(e_i,\cdot)=\nabla F_0(e_i,\cdot)\ \forall\ i=1,2,\dots,d])>1-\e;\\
\end{align*}
\sms {\rm (b)} $\exists$ a $\mathcal R_T$-holonomy $R$ with $\mathcal D(R),\ \mathcal I(R)\subset A,\ \ \&\ \mu(\mathcal D(R))\ge r\mu(A)$ so that
$F(R(x))-F(x)=\sigma\ \forall\ x\in\mathcal D(R)$.
\endproclaim
\demo{Proof}

Since $F_0:X\to\Bbb Z$ is bounded,\ $\exists\ K\in\Bbb N,\ A_1,A_2,\dots,A_K\in\B\cap A$ so that
$A=\bigcupdot_{k=1}^KA_k$ so that $F_0$ is constant on each $A_k$.
\

Consider the measurable partition $\a:=\{A_1,A_2,\dots,A_K,X\setminus A\}$ of $X$.
\

Fix $0<\d\ll\e$ and in particular $\d<\tfrac14$. By the ergodic theorem,
$$\exists\ \frak n\ge N_{\mathcal T_0}\ \text{so that}\ \ m(\frak a_\frak n)>1-\d^2$$ where
$$\frak a_\frak n:=\bigcap_{a\in \a,\ n\ge \frak n}\left[\frac1{|\Si_n|}\sum_{k\in \Si_n}1_a\circ T_k=m(a)(1\pm\d)\right].$$
\

Now let $\mathcal T$ be a Rokhlin tower with $N_\mathcal T\ge \frac{2\frak n}{\d}$ and
$\e_\mathcal T<\d$.
\

Let $\widehat{\a}=\a\vee\{\frak a_\frak n,\frak a_\frak n^c\}\vee\{B_{\mathcal P_0},B_{\mathcal P_0}^c\}$ and let
$\beta $ be the $\widehat{\a}$-purification of $B_\mathcal T$ with $\mathcal P$ the   $\widehat{\a}$-purification of $\mathcal T$.

 Define
$$\beta _{\text\tiny\smiley}:=\{b\in \beta :\ \#\{k\in \Si_{N_\mathcal T}:\ T_kb\subset\frak a_\frak n\}>(1-\d)|\Si_{N_\mathcal T}|\}\ \&\ U_{\text\tiny\smiley}:=\bigcupdot_{b\in\beta _{\text\tiny\smiley}}b.$$
By the Chebyshev-Fubini theorem,
$$m( U_{\text\tiny\smiley})>(1-\tfrac{\d}{1-\d})\cdot m(B_\mathcal T)>\tfrac12\cdot m(B_\mathcal T).$$

To obtain $F$ satisfying the essential value condition, we make two  changes to $F_0$. The first preparatory change is to ensure that $F$ will be $\mathcal T$-internal and
$\mathcal S$-incremental. Let
 $$Q:=\Si_{N_\mathcal T}\cap\Si(\bdy\Si_{N_{\mathcal T}},N_{\mathcal T_0})\sqcup\Si(\bdy\Si_{\lfl\frac{N_{\mathcal T}}2\rfl},N_{\mathcal T_0})$$
and define
$$F_1(x):=\begin{cases} & 0\ \ \ \ \ \ x\in T_kb\ \text{where}\ b\in \beta, k\in Q \ \&\ \exists\ k_0\in Q\cap\Sigma(k,N_{\mathcal T_0}),\
T_{k_0}b\subset B_{\mathcal T_0},\\ &
F_0(x)\ \ \ \ \ \ \text{else,}\end{cases}$$
then
$$m([F_0\ne F_1])\le \# Qm(B_\mathcal T).$$
Now,  for some constant $C=C_d>0$,
$$\# Q\le CN_{\mathcal T_0}N_\mathcal T^{d-1}\ \&\ m(B_\mathcal T)\le \frac{C}{N_\mathcal T^d}$$
whence
$$m([F_0\ne F_1])\le 
 \frac{C^2N_{\mathcal T_0}}{N_{\mathcal T}}.$$
The change was made on full $\mathcal T_0$ subtowers in the purification and so
 $F_1|_{T_QB_{\mathcal T}}\equiv 0$ since $F_0$ is $\mathcal T_0$-internal. In particular, $F_1$ is $\mathcal T$-internal.
Moreover, $F_1$ is $\mathcal S$-incremental again because the change was made on full $\mathcal T_0$ subtowers in the purification and on the remaining $\mathcal T_0$ subtowers,
$F_0$ is $\mathcal S$-incremental.
\

We can now define $F:X\to\Bbb G$ by
$$F(x)=\begin{cases}& F_1(x)+\sigma\ \ \ \ \ \ \ x\in T_{\Si_{\lfl\frac{N_{\mathcal T}}2\rfl}}B_\mathcal T\\ & F_1(x)\ \ \ \ \ \ \ \text{else}.\end{cases}$$
The changes made were small enough for (a).
It follows that $F$ is $\mathcal T_0$-internal and $\mathcal S$-incremental.

\demo{ Proof {\rm that $F$ satisfies (b)}}
\

 For each  $b\in\beta_{\smiley},\ j\in \{1,2,\dots,K\}$, let
 $$K_{b,j}:=\{k\in \Sigma_{N_{\mathcal T}}:\ T_kb\subset A_j\}.$$
 \
 
 By construction,  for each  $b\in\beta\ \&\ j\in \{1,2,\dots,K\}$,
 $$\# K_{b,j}\cap\Sigma_{\frac{N_{\mathcal T}}2}\ll K_{b,j}\cap \Sigma_{N_{\mathcal T}}
 \setminus\Sigma_{\frac{N_{\mathcal T}}2}$$ and there is an injection 
 $$\frak k_{b,j}:K_{b,j}\cap\Sigma_{\frac{N_{\mathcal T}}2}\to K_{b,j}\cap \Sigma_{N_{\mathcal T}}
 \setminus\Sigma_{\frac{N_{\mathcal T}}2}$$
 Define the $\mathcal R_T$-holonomy by
 $$R(x)=T_{\frak k_{b,j}(j)-j}(x)\ \ x\in T_jb,\ b\in\b,\ j\in K_{b,j}\cap\Sigma_{\frac{N_{\mathcal T}}2}\},$$
 then $F(R(x))-F(x)=\s$ and 
                            
 \begin{align*}\mu(\mathcal D(R))&\ge
 \sum_{b\in\beta_{\smiley}}\sum_{j=1}^K\sum_{k\in K_{b,j}\cap\Sigma_{\frac{N_{\mathcal T}}2}}\mu(A\cap T_kb)\\ &=
 \int_{U_{\smiley}}\sum_{k\in \Sigma_{\frac{N_{\mathcal T}}2}}1_A\circ T_kdm\\ &>
 (1-\d)|\Sigma_{\frac{N_{\mathcal T}}2}|m(U_{\smiley})m(A)\\ &\ge
 (1-2\d)|\Sigma_{\frac{N_{\mathcal T}}2}|m(B_\mathcal T)m(A)\\ &\ge
(1-2\d)(1-\d)\tfrac{|\Sigma_{\frac{N_{\mathcal T}}2}|}{|\Sigma_{N_{\mathcal T}}|}m(A)\\ &\ge 
\tfrac{(1-2\d)(1-\d)}{3\cdot 2^{d}}m(A)>rm(A). \ \ \CheckedBox
\end{align*}
To finish the proof of theorem 2,  fix
\bul a countable, dense collection $\mathcal A\subset\B$;
\bul a finite, symmetric $\mathcal S\subset\Bbb G$ so that $\overline{\<\mathcal S\>}=\Bbb G$.
Write down a sequence
$$((\sigma_n,A_n):\ n\in\Bbb N)\in(\mathcal S\x\mathcal A)^{\Bbb N}$$ so that
for each $(\sigma,A)\in \mathcal S\x\mathcal A$,
$$\#\,\{n\in\Bbb N:\ (\sigma_n,A_n)=(\sigma,A)\}=\infty.$$
Now apply the Inductive Lemma recursively with  $r=\tfrac1{2^{d+4}}$ and $0<\e_{n+1}<\e_n<\tfrac1{2^n}$ to obtain a sequence of
 $\mathcal S$-incremental $F_n:X\to\Bbb G$ so that for each $n\ge 1$,

{\small\begin{align*}&\tag{a$_n$}
 \mu([\nabla F_{n}(e_i,\cdot)\ne\nabla F_{n-1}(e_i,\cdot)\ \forall\ i=1,2,\dots,d])<\e_n;\\
 &\tag{b$_n$}\mu(\{x\in A_n:\ \exists\ u\in\Bbb Z^d,\ T_ux\in A_n\ \&\ \nabla F_n(u,x)=\sigma_n\})\ge rm(A_n);\\ &\tag{c$_n$}
 G:X\to\Bbb G\ \&\ \mu([\nabla G(e_i,\cdot)\ne\nabla F_{n}(e_i,\cdot)\ \forall\ i=1,2,\dots,d])<\e_{n+1}\\ &\ \ \ \ \ \ \ \ \ \ \ \ \ \  \Lra\ G\ \text{satisfies (b$_k$)}\ \forall\ 1\le k\le n-1.
\end{align*}}
It follows that for $m$-a.e. $x\in X,\ k\in\Bbb Z^d,\ \exists\ \lim_{n\to\infty}F_n(k,x)=:F(k,x)$ and $F:\Bbb Z^d\x X\to\Bbb G$ is a $\mathcal S$-incremental cocycle
satisfying   $$\text{\tt EVC}_T(\sigma+U,C,A )\ \forall\ A\in \mathcal A,\ \sigma\in\mathcal S,\ U\in\frak T_0.$$
By the Ergodicity Proposition, the skew product action is ergodic.\ \Checkedbox
\section*{\S3 Infinitely generated groups}
\

We do not know if theorem 1 holds for all infinite, countable groups.
In this section we prove versions  for certain examples of infinitely generated groups.

\subsection*{{\sc Locally Finite Groups}}
\

Say that $\G$ is a  {\it normally, locally finite group} if it is the increasing union of finite normal subgroups
$$G_1 \triangleleft G_2 \triangleleft G_3 \triangleleft\dots\uparrow\G.$$
\proclaim{Theorem 3} Let $\G$ be a normally,   locally finite group, let  $(X,d)$ be a perfect Polish space  and let
$T:\G\to\text{\tt Homeo}\,(X)$ be a
{\tt tt} action.
\

For any separable Banach space, $\exists$ a continuous cocycle $h:\G\x X\to\Bbb B$ so that the skew product action
$T^{(h)}:\G\to\text{\tt Homeo}\,(X\x\Bbb B)$ is {\tt tt}.\endproclaim\demo{Proof}\ Let
$x_0\in X$ be  a properly recurrent point, i.e. $x_0\in\overline{T_{\G\setminus\{e\}}(x_0)}$.

\

We claim that because of the finiteness of the groups $G_n$ and the perfectness of $X$,
\Par \ \ Given  $U\in\mathcal T_{x_0},\ k\ge 1,\ \exists\  N >k\ \&\ \g_U\in G_N\setminus G_k$ so that $T_{\g_U}(x_0) \in U\setminus\{x_0\}$.
\

We'll need in  addition, the following lemma:
\proclaim{{\Large\dsrailways}}\ \ Let $F\in C(X,\Bbb B),\ s\in\Bbb B,\ \D>0$, then
  $\exists\ f\in C(X,\Bbb B)$ so that
\begin{align*}&\tag{i}f\circ T_g\equiv f\ \forall\ g\in G_{N-1};\\
 &\tag{ii}\bdy (F+f)\  \text{satisfies}\
\text{\tt EVC}(x_0,U,s,\D,\g_U).
\end{align*}\endproclaim
\demo{Proof of \dsrailways}\ \ By possibly shrinking $U$ (and suitably adjusting $\g_U\ \&\ N$) we can ensure that

$$\sup_{y\in U}\,\|F(y)-F(x_0)\|_\Bbb B<\D.$$
\

Fix $\mathcal E>0$ so that the sets
 $$\{T_gB(x_0,2\mathcal E):\ g\in G_N\}$$ are disjoint. Here, for $x\in X,\ r>0$,
 $B(x,r):=\{y\in X:\ d(y,x)\le r\}$ is the closed $r$-ball centered at $x$.
 \

For $x\in X,\ g\in G_N$, let
 \begin{align*}
    w(x):=\left(1-\frac{d(B(T_{\g_U}(x_0),\mathcal E),x)}{\mathcal E}\right)_+.
 \end{align*}
 where $d(B,x):=\inf_{y\in B}d(y,x)$ and $a_+:=\max\{a,0\}$ for $a\in\Bbb R$.
It follows that
$$w\in  C(X,[0,1]),\ {w}|_{B(T_{\g_U}(x_0),\mathcal E)}\equiv 1\ \&\ {w}|_{B(T_{\g_U}(x_0),2\mathcal E)^c}\equiv 0.$$
Now define $f:X\to\Bbb B$ by
$$g(x)=\begin{cases} &\ sG(T^{-1}_{g}x)\ \ \ \ \ \ \ \ x\in T_gB(T_{\g_U}(x_0),2\mathcal E),\ g\in G_{N-1},\\ &
  0\ \ \ \ \ \ \ \ x\in X\setminus\bigcup_{g\in G_{N-1}}T_gB(x_0,2\mathcal E). \end{cases}$$

This $f\in  C_B(X,\Bbb B)$ is as required $\because\ \ \g_UG_{N-1}=G_{N-1}\g_U$.\ \ \Checkedbox\ \dsrailways
\

 \demo{Proof of theorem 3}
\

Fix
\bul a decreasing sequence $(U_n:\ n\ge 1),\ U_n\in\frak T_{x_0}$ so that $\forall\ V\in\frak T_{x_0}\ \exists\ N_V$ so that
$U_n\subset V\ \forall\ n\ge N_V$;
\bul $(s_n:\ n\ge 1) \in \Bbb B^\Bbb N$ dense, and taking each value i.o.;
\bul $\D_n\downarrow 0$.
\

Using \ \dsrailways   \ iteratively, construct:
\sms\ $\kappa_n<\kappa_{n+1}\to\infty$ and $\g_{{U_n}}\in G_{\kappa_n}\setminus G_{\kappa_{n-1}}$ so that
$T_{\g_{{U_n}}}(x_0)\in U_n$;
\sms\ $f_n\in C(X,\Bbb B)$ so that
\begin{align*}&\tag{a}f_n\circ T_g\equiv f_n\ \forall\ g\in G_{\kappa_{n-1}};\\
 &\tag{ii}\bdy \sum_{k=1}^nf_k\  \text{satisfies}\
\text{\tt EVC}(x_0,U_n,s,\D_n,\g_{U_n}).
\end{align*}
For each $\g\in\G,\ x\in X$, the sum
$$h(\g,x):=\sum_{n=1}^\infty (f_n(T_\g(x))-f_n(x))$$
converges (only finitely many elements being non-zero);
$h:\G\x X\to\Bbb B$ is a continuous cocycle.
\

Now for $\g\in G_{\kappa_{N}}$,
$h(\g,x)=\sum_{n=1}^N (f_n(T_\g(x))-f_n(x))$, whence
 $h$ satisfies
$\text{\tt EVC}\,(x_0,U_n,s_n,0,\g_{{U_n}})\ \forall\ n\ge 1$.

\

By proposition 1,
  $E(h,x_0)=\Bbb B\ \&$ hence $T^{(h)}$ is {\tt tt}.\ \  \Checkedbox\ 
  \subsection*{\sc Actions of $\Bbb Z^\infty$}
  \

  Here we consider
  $$\Bbb Z^\infty:=\{g=(g_1,g_2,\dots)\in\Bbb Z^\Bbb N:\ g_n=0\ \forall\ n\ \text{\tt large}\}$$
  with coordinatewise addition. The multiplicative group $\Bbb Q_+$ is isomorphic with $\Bbb Z^\infty$ by prime factorization.

  \

  For $F\subset\Bbb N$ finite,
 define $\pi_F:\Bbb Z^\infty\to\Bbb Z^F$ by $\pi_F(\g):=\g|_F$.

\proclaim{Proposition 4}
\

Let $(X,d)$ be a perfect polish space,  let $T:\Bbb Z^\infty\to\text{\tt Homeo}\,(X)$ be a free,
{\tt tt} action, and let $\Bbb B$ be a separable Banach space.
 \

 There is a cocycle $h:\Bbb Z^\infty\x X\to\Bbb B$ so that the skew product action $T^{(h)}:\Bbb Z^\infty\x (X\x\Bbb B)\to\text{\tt Homeo}\,(X\x\Bbb B)$ is
 {\tt tt}.\endproclaim\demo{Proof} There are two cases covering the theorem.
 \f{\bf Case 1} There is a $T$-transitive  point $x_0\in X$ and $N\ge 1,\ \g_k\in\Bbb Z^{[1,N]}$ so that
 $T_{\g_k}(x_0)\xrightarrow\ x_0$.
 \demo{Proof of the proposition in case 1} Since $\Bbb Z^{[1,N]}$ is finitely generated, by the  proof of theorem 1,
 there is a cocycle $\eta:\Bbb Z^{[1,N]}\x X\to\Bbb B$ so that
 $$\overline{\{T^{(\eta)}_g(x_0,0):\ g\in\Bbb Z^{[1,N]}\}}\supset \{x_0\}\x\Bbb B.$$
 Define $h:\Bbb Z^\infty\x X\to\Bbb B$ by
 $$h(\g,x):=\eta(\pi_{[1,N]}(\g),x),$$ then $h:\Bbb Z^\infty\x X\to\Bbb B$ is a cocycle and
 $$\overline{\{T^{(h)}_\g(x_0,0):\ \g\in\Bbb Z^\infty\}}=\overline{\{T^{(\eta)}_g(x_0,0):\ g\in\Bbb Z^{[1,N]}\}}\supset \{x_0\}\x\Bbb B.$$
 By transitivity of $x_0$ for $T$, we have
 $$\overline{\{T^{(h)}_\g(x_0,0):\ \g\in\Bbb Z^\infty\}}=X\x\Bbb B.\ \ \CheckedBox$$
 If Case 1 fails, then we are in
 \f{\bf Case 2} There is a $T$-transitive  point $x_0\in X$ and $\g_k\in\Bbb Z^\infty$ so that $T_{\g_k}(x_0)\xrightarrow\ x_0$ and
 $\forall\ k,\ K\ge 1,\ exists\ L>K$ so that $(\g_k)_L\ne 0$.
  \demo{Proof of the proposition in case 2} \ \ We prove
 \proclaim{\smiley}\  $\exists$ a homomorphism  $H:\Bbb Z^\infty\to\Bbb B$ so that the direct product action $T\x H:\g\to\text{\tt Homeo}\,(X\x\Bbb B)$ is
 {\tt tt}.\endproclaim

 \

\demo{Proof of \smiley}

In the absence of case 1, we have that
\begin{align*}\tag*{\Industry}&\forall\ V\in\frak T_{x_0},\ s\in\Bbb B\ \&\  n_0\ge 1,\ \exists\\ &
 n>n_0\ \&\ \g\in\Bbb Z^{[1,n]}\ \text{with}\ \g_n\ne 0\ \&\ T_\g(x_0)\in V.
\end{align*}
Given $n,\ \g\ \&\ s$ as above, let  $h_{n,\g,s}:\Bbb Z^\infty\to\Bbb B$ be the homomorphism given by
$$h_{n,\g,s}(g)=\pi_{n,\g}(g_n)\cdot s:=\tfrac{g}{\g_n}\cdot s.$$

\

Fix
\bul a decreasing sequence $(U_n:\ n\ge 1),\ U_n\in\frak T_{x_0}$ so that $\forall\ V\in\frak T_{x_0}\ \exists\ N_V$ so that
$U_n\subset V\ \forall\ n\ge N_V$;
\bul $(s_n:\ n\ge 1) \in \Bbb B^\Bbb N$ dense, and taking each value i.o. .
\

Using \ \Industry,\  \ iteratively construct:
\sms\ $\kappa_n<\kappa_{n+1}\to\infty$ and $\g^{(n)}\in \Bbb Z^{[1,\kappa_n]}),\ \g^{(n)}_{\kappa_n}\ne e_{G_{\kappa_n}}$ so that
$T_{\g^{(n)}}\in U_n$.
\sms\ We have that
for each $\g\in\Bbb Z^\infty$, the sum
$$H(\g):=\sum_{n=1}^\infty h_{\kappa_n,\g^{(n)},s_n}(\g)$$
converges (only finitely many elements being non-zero);
$H:\Bbb Z^\infty\to\Bbb B$ is a homomorphism.
\

Considering $H:\Bbb Z^\infty\x X\to\Bbb B$ as a (constant) cocycle, we have that $H$ satisfies
$\text{\tt EVC}\,(x_0,U_n,s_n,0,\g^{(n)})\ \forall\ n\ge 1$.

\

By proposition 1,  we have that
  $E(H,x_0)=\Bbb B\ \&$ hence $T^{(H)}=T\x H$ is {\tt tt}.\ \ \ \ \Checkedbox\ \smiley

  \section*{\S4 H\"older continuous cocycles for $\Bbb Z^d$ shifts}
  Let $S$ be a finite set, let $d\ge 1\ \&$ let $X=S^{\Bbb Z^d}$.
\

 The function $f:X\to\Bbb R^D$ is {\it H\"older continuous} if for some $\th\in (0,1)\ \&\ M>0$,
$$\|f(x)-f(y)\|_2\le M\th^{t(x,y)}\ \forall\ x,\ y\in X$$
where   $t(x,y):=\min\,\{\|n\|_1:\ n\in\Bbb Z^d,\ x_n\ne y_n\}$. Here and throughout, for
$N\ge 1\ \&\ p>0$,
$\|(x_1,\dots,x_N))\|_p:=(\sum_{k=1}^N|x_k|^p)^\frac1p$.
\

A H\"older continuous function taking finitely many values is aka a {\it block function}.
\

A cocycle $F:\Bbb Z^d\x X\to\Bbb R^D$ is called H\"older continuous   if $x\mapsto F(n,x)$
is H\"older continuous $\forall\ n\in\Bbb Z^d$.
\

\subsection*{Example:\ Random walks}
\

\ Let $D\ge 1$ and let $(X,T)=(S^\Bbb Z,\text{\tt Shift})$ where $S$ is  finite set, large enough so that
     $\exists$ $\phi:S\to\Bbb R^D$ with
     $\overline{\text{\tt Semigroup}(\phi(S))}=\Bbb R^D$.
     \

     Define $\varphi:X\to\Bbb R^D$ by
$\varphi(x):=\phi(x_0)$. It is not hard to see that the skew product  $T_\varphi$ is {\tt tt}.
\

If $\mu\in\mathcal P(S)$ satisfies $\sum_{s\in S}\mu(s)\varphi(s)=0$ then $(X\x\Bbb R^D,\mu^\Bbb Z\x m_{\Bbb R^D},T_\varphi)$ is a
measure preserving transformation and is ergodic if $D=1, 2$ (see \cite{Harris-Robbins} ).
 For $D\ge 3,\ T_\varphi$ is dissipative
 by the local limit theorem (see e.g. \cite{Breiman}) whence not ergodic.

These random walk constructions  have no ergodic analogues for higher dimensional actions.
The reason is basically
\proclaim{Schmidt's theorem}
Let $(X,T)=(S^{\Bbb Z^d},\text{\tt Shift})$ where $d\ge 2$ and $S$ is a finite set and let 
$F:\Bbb Z^d\x X\to\Bbb R^k$ be a H\"older continuous cocycle, then
$$F(n,x)=g(T_nx)-g(x)+H(n)$$
where $g:X\to\Bbb R^k$ is H\"older continuous and $H:\Bbb Z^2\to\Bbb R^k$ is  a homomorphism.\endproclaim
This can be deduced from the more general theorem 3.2 in \cite{Schmidt-HOcohom}
which is a symbolic version of  a similar result for multidimensional Anosov actions
(theorem 2.9 in \cite{Katok-Spatzier}).
\

\proclaim{Corollary}\ \ Let $d\ge 2\ \&
 (X,T)=(S^{\Bbb Z^d},\text{\tt Shift})$ where $S$ is a finite set.
For $d'\ge d,$ there is no H\"older continuous cocycle $F:\Bbb Z^d\x X\to\Bbb R^{d'}$ with the skew product
$(X\x\Bbb R^{d'},T^{(F)})$ topologically transitive.\endproclaim
\demo{Proof}  \ \ Let $F:\Bbb Z^d\x X\to\Bbb R^{d'}$ be a H\"older continuous cocycle.  By Schmidt's theorem
$$F(n,x)=g(T_nx)-g(x)+H(n)$$
where $g:X\to\Bbb R^{d'}$ is H\"older continuous and $H:\Bbb Z^d\to\Bbb R^{d'}$ is  a homomorphism, whence
the skew product action $T^{(F)}$ is continuously conjugate to the product action
 $T\x H$ where $H_n(z):=z+H(n)$.
 \

 Topological transitivity is impossible since either $\dim\text{Span}\,H(\Bbb Z^d)<d'$
or $H(\Bbb Z^d)$ is discrete. \ \Checkedbox
 \

On the other hand, 
\Smi\ There is  an homomorphism $H:\Bbb Z^d\to\Bbb R^{d-1}$ whose  product action
$(X,m\x m_{\Bbb R^{d-1}},T\x H)$ is
ergodic.
\

\demo{Proof of {\large\smiley}}\ \ Let $d\ge 1$, let $(X,T)=(\{0,1\}^{\Bbb Z^d},\text{\tt Shift})$  and let
$m\in\mathcal P(\{0,1\}^{\Bbb Z^d})$ be symmetric product measure.
\

Define $H:\Bbb Z^d\to\Bbb R^{d-1}$ by
$$H((n_1,n_2,\dots,n_d)):=\sum_{j=1}^{d-1}n_je^{(d-1)}_j+n_d\vec{\a}$$
with $\vec{\a}:=\sum_{k=1}^{d-1}\alpha_ke^{(d-1)}_k$
where $1,\a_1,\a_2,\dots,\a_{d-1}$ are linearly independent over $\Bbb Q$ so that $(\Bbb T^{d-1}, m_{\Bbb T^{d-1}},R_{\vec{\a}})$
is ergodic ($R_{\vec{\a}}(x):=x+\vec{\a})$).
\

The action $(X,T,m)$ is strongly mixing, whence {\tt mildly mixing}. By \cite{S-W},
the product $\Bbb Z^d$ action $(X\x\Bbb R^{d-1},T\x H,m\x m_{\Bbb R^{d-1}})$ is ergodic if and only if
the action $(\Bbb R^{d-1},H,m_{\Bbb R^{d-1}})$ is conservative, ergodic.
\

We have that $(\Bbb R^{d-1},m_{\Bbb R^{d-1}},H)\cong (\Bbb T^{d-1}\x\Bbb Z^{d-1},m_{\Bbb T^{d-1}}\x\#,J)$
where $\cong$ denotes measure theoretic isomorphism of actions and
$$J_{e_k^{(d)}}(x,z):=\begin{cases}& (x,z+e_k^{(d-1)})\ \ \ \ \ \ \ \ \ 1\le k\le d-1,\\ &
                (x+\vec{\a},z+\tau(x))\ \ \ \ \ \ \ \ \ k=d
                \end{cases}$$
                with

                $$\tau(x):=(\lfl x_1+\a_1\rfl,\lfl x_2+\a_2\rfl,\dots,\lfl x_{d-1}+\a_{d-1}\rfl).$$
               If $F:\Bbb T^{d-1}\x\Bbb Z^{d-1}\to\Bbb C$ is a measurable, $J$-invariant function  ($F\circ J_{n}=F\  \forall\ n\in\Bbb Z^d)$) then
$$F(x,n)=F\circ J_{n,0}(x,0)=F(x,0)\ \&\ F(x,0)=F\circ J_{e_d}(x,0)=F(R_{\vec{\a}}(x),0).$$

Thus, by ergodicity of $R_{\vec{\a}}$, $F$ is a.e. constant and $J$ is ergodic. Conservativity follows as the underlying measure
space is non-atomic.\ \ \Checkedbox
\subsection*{Concluding Remarks}\ We have been unable to establish versions of theorem 1   
for actions of $\Bbb Q$ or $\Si(\Bbb N)$  (the group of finite permutations of $\Bbb N$). 
We have also been unable to prove theorem 2 for cocycles with values in an arbitrary, countable amenable group (established for $\Bbb Z$-actions
in \cite{AW-Herman}).

\

On the other hand,
theorem 2 can be generalized to free actions of countable
amenable groups.  The appropriate inductive lemma  is also established using the ergodic theorem and
 an appropriate Rokhlin lemma   for countable
amenable group actions as in
 \cite{O-W}.


\begin{thebibliography}{{Bes}37}

\bibitem[ALV98]{A-L-V}
Jon Aaronson, Mariusz Lemanczyk, and Dalibor Volny.
\newblock A cut salad of cocycles.
\newblock {\em Fund. Math.}, 157(2-3):99--119, 1998.
\newblock Dedicated to the memory of Wieslaw Szlenk.

\bibitem[AW04]{AW-Herman}
Jon Aaronson and Benjamin Weiss.
\newblock On {H}erman's theorem for ergodic, amenable group extensions of
  endomorphisms.
\newblock {\em Ergodic Theory Dynam. Systems}, 24(5):1283--1293, 2004.

\bibitem[{Bes}37]{Bes37}
A.~S. {Besicovitch}.
\newblock {A problem on topological transformation of the plane.}
\newblock {\em {Fundam. Math.}}, 28:61--65, 1937.

\bibitem[Bre68]{Breiman}
Leo Breiman.
\newblock {\em Probability}.
\newblock Addison-Wesley Publishing Company, Reading, Mass.-London-Don Mills,
  Ont., 1968.

\bibitem[CFW81]{C-F-W}
A.~Connes, J.~Feldman, and B.~Weiss.
\newblock An amenable equivalence relation is generated by a single
  transformation.
\newblock {\em Ergodic Theory Dynamical Systems}, 1(4):431--450 (1982), 1981.

\bibitem[Con73]{Conze}
J.~P. Conze.
\newblock Entropie d'un groupe ab\'elien de transformations.
\newblock {\em Z. Wahrscheinlichkeitstheorie und Verw. Gebiete}, 25:11--30,
  1972/73.

\bibitem[FM77]{F-M}
Jacob Feldman and Calvin~C. Moore.
\newblock Ergodic equivalence relations, cohomology, and von {N}eumann
  algebras. {I}.
\newblock {\em Trans. Amer. Math. Soc.}, 234(2):289--324, 1977.

\bibitem[Her79]{Herman}
M.R. Herman.
\newblock Construction de diff\'eomorphismes ergodiques.
\newblock {\em Unpublished Manuscript}, 1979.

\bibitem[HR53]{Harris-Robbins}
T.~E. Harris and Herbert Robbins.
\newblock Ergodic theory of {M}arkov chains admitting an infinite invariant
  measure.
\newblock {\em Proc. Nat. Acad. Sci. U. S. A.}, 39:860--864, 1953.

\bibitem[KS94]{Katok-Spatzier}
Anatole Katok and Ralf~J. Spatzier.
\newblock First cohomology of {A}nosov actions of higher rank abelian groups
  and applications to rigidity.
\newblock {\em Inst. Hautes \'Etudes Sci. Publ. Math.}, (79):131--156, 1994.

\bibitem[KW72]{K-W}
Yitzhak Katznelson and Benjamin Weiss.
\newblock Commuting measure-preserving transformations.
\newblock {\em Israel J. Math.}, 12:161--173, 1972.

\bibitem[OW87]{O-W}
Donald~S. Ornstein and Benjamin Weiss.
\newblock Entropy and isomorphism theorems for actions of amenable groups.
\newblock {\em J. Analyse Math.}, 48:1--141, 1987.

\bibitem[Sch77]{S}
Klaus Schmidt.
\newblock {\em Cocycles on ergodic transformation groups}.
\newblock Macmillan Company of India, Ltd., Delhi, 1977.
\newblock Macmillan Lectures in Mathematics, Vol. 1.

\bibitem[Sch95]{Schmidt-HOcohom}
Klaus Schmidt.
\newblock The cohomology of higher-dimensional shifts of finite type.
\newblock {\em Pacific J. Math.}, 170(1):237--269, 1995.

\bibitem[Sid73]{Sidorov}
E.~A. Sidorov.
\newblock Topologically transitive cylindrical cascades.
\newblock {\em Mat. Zametki}, 14:441--452, 1973.
\newblock English translation: Math. Notes 14 (1973), 810–816 (1974).

\bibitem[Sni30]{Sni30}
L.~G. Snirelman.
\newblock {Example of a transformation of the plane}.
\newblock {\em {Izv. Donsk. Politekhn. Inst.}}, 4:64--74, 1930.

\bibitem[SW82]{S-W}
Klaus Schmidt and Peter Walters.
\newblock Mildly mixing actions of locally compact groups.
\newblock {\em Proc. London Math. Soc. (3)}, 45(3):506--518, 1982.

\end{thebibliography}
\end{document}